\documentclass[12pt]{article}
\usepackage{amsmath, amssymb, amsthm, amsfonts, enumerate}
\usepackage{graphicx}

\newtheorem{theorem}{Theorem}

\def\al{\alpha}
\def\be{\beta}
\def\ga{\gamma}
\def\de{\delta}

\def\si{\sigma}

\def\noi{\noindent}

\def\proof{\noi {\bf Proof.} }
\def\gcd{{\rm gcd} \, }

\begin{document}

\title{Tiling of regular polygons with similar right triangles}

\author{Mikl\'os Laczkovich and Ivan Vasenov}

\footnotetext[1]{{\bf Keywords:} Tilings with right triangles, regular polygons}
\footnotetext[2]{{\bf MR subject classification:} 52C20}
\footnotetext[3]{The first author was supported by the Hungarian National
Foundation for Scientific Research, Grant No. K124749.}

\maketitle

\begin{abstract} 
We prove that for every $N\ne 4$ there is only one right triangle that tiles the regular $N$-gon.
\end{abstract}

We say that a triangle $T$ tiles a polygon $A$, if $A$ can be dissected into
finitely many nonoverlapping triangles similar to $T$.

Let $R_N$ denote the regular $N$-gon. Connecting the center of $R_N$ with the
vertices of $R_N$ we obtain a dissection of $R_N$ into $N$ congruent isosceles
triangles. Bisecting each of these triangles into two right angled triangles, we
get a dissection of $R_N$ into $2N$ congruent right triangles with acute angles
$\pi /N$ and $(\pi /2)-(\pi /N)$. In particular, we find that the
right triangle with acute angles $\pi /N$ and $(\pi /2)-(\pi /N)$ tiles $R_N$.

In this note we are concerned with the following question: 
are there other right triangles that tile $R_N$? If $N=4$, then there
are infinitely many such triangles. Indeed, the right triangle with legs $1$ and
$1/n$ tiles the unit square for every positive integer $n$. In \cite{L1} and
\cite{Sz} it is proved that a right triangle with acute angle $\al$ tiles the
square if and only if $\tan \al$ is a totally positive algebraic number; that
is, if every real conjugate of $\tan \al$ is positive. 

It is also known that if $N\ge 25$ and $N\ne 30,42$, then the only right
triangle that tiles $R_N$ has acute angles $\pi /N$ and $(\pi /2)-(\pi /N)$
(see \cite[Theorem 1.1]{L2}). The same was proved for $N=5$ in \cite{SH}.
We prove that this is true for every $N\ne 4$.

\begin{theorem}\label{t1}
For every $N\ne 4$ there is only one right triangle that tiles the regular
$N$-gon.
\end{theorem}

\proof
Let $T$ be a right triangle with acute angles $\al ,\be$. We prove that if $T$
tiles $R_N$, where $N\ne 4$, then one of $\al$ and $\be$ equals $\pi /N$.

Let $\de _N$ denote the angle of $R_N$; that is, let $\de _N =(N-2)\pi /N$.
Suppose that $T$ tiles $R_N$, where $N\ne 4$, and fix such a tiling. Let $V$
be a vertex of any of the tiles. We say that $p\al +q\be +r\ga =\si$ is the
equation at $V$ if the number of triangles having $V$ as a vertex and having
angle $\al$ (resp. $\be$ or $\ga =\pi /2$) at $V$ equals $p$ (resp. $q$ or
$r$).  Here $\si =\de _N$ if $V$ is one of the vertices of $R_N$, and $\si =\pi$
or $2\pi$ otherwise.

Let the equation at any of the vertices of $R_N$ be $p\al +q\be +r\ga =\de _N$.
Since $N\ne 4$, $\de _N$ is not an integer multiple of $\pi /2$, and thus
$p\ne q$. We obtain $(p-q)\al +(q+r)\pi /2=\de _N =(N-2)\pi /N$, showing that
$\al$ is a rational multiple of $\pi$. Then so is $\be =(\pi /2)-\al$.

Suppose $N\ne 6$. By Theorem 1.2 of \cite{L2}, each angle of $R_N$ is packed
with at most two tiles. Suppose $N=3$. Then we have either
$\al =\de _3 = \pi /3$, or $2\al =\de _3 = \pi /3$, hence $\al =\pi /6$ and
$\be =\pi /3$. This proves the statement for $N=3$.

Therefore, we may assume $N\ge 5$. Then $\de _N >\pi /2$, and there must be
exactly two tiles at each vertex of $R_N$. Then we have either $2\al =\de _N =
\pi - (2\pi /N)$ or $(\pi /2)+\al =\pi -(2\pi /N)$.

If $2\al =\de _N$, then we get $\al =(\pi /2)-(\pi /N)$ and $\be =\pi /N$,
and we are done.

Suppose $(\pi /2)+\al =\pi -(2\pi /N)$. Then $\al =(\pi /2)-(2\pi /N)=
(N-4)\pi /(2N)$ and $\be = 2\pi /N$. By Lemma 1.10 of \cite{L2}, we have
\begin{equation}\label{e1}
\left\{ \frac{k(N-4)}{2N}\right\} + \left\{ \frac{2k}{N}\right\} +
\left\{ \frac{k}{2}\right\} =1
\end{equation}
whenever $\gcd(k,2N)=1$ and $\{ k/N\} <1/2$. Let
$$k=
\begin{cases}
(N/2)-1 &\text{if $N \equiv 0$ (mod $4$)},\\
N+(N-1)/2 &\text{if $N \equiv 1$ (mod $4$)},\\
(N/2)-2 &\text{if $N \equiv 2$ (mod $4$)},\\
(N-1)/2 &\text{if $N\equiv 3$ (mod $4$)}.\\
\end{cases}
$$
Then $\gcd (k,2N)=1$ and $\{ k/N\} <1/2$. We have $\{ k/2\} =1/2$, since $k$ is
odd. Then \eqref{e1} gives $\{ 2k/N\} \le 1/2$. However, it is easy to check
that $\{ 2k/N\} >1/2$ for every $N\ge 5$, $N\ne 6$. Thus the case $\al =
(\pi /2)-(\pi /N)$ is impossible. This proves the statement for every
$N\ne 4,6$.

Finally, suppose $N=6$. In the following argument we assume $\al \le \be$. Put
$b=2\be /\pi$; then we have $1/2\le b<1$. If the equation at a vertex of $R_6$
is $p_0 \al +q_0 \be +r_0 \ga =\de _6 =2\pi /3$, then we have
\begin{equation}\label{e2}
p_0 (1-b)+q_0 b+r_0 =4/3.
\end{equation}
Suppose $p_0 \le q_0 $. Then $(q_0 -p_0 )b+p_0 +r_0 =4/3$ gives $sb=4/3$ or
$sb=1/3$, where $s=q_0 -p_0$ is a nonnegative integer. Since
$1/2 \le b< 1$, we have $s=2$, $sb=4/3$, $b=2/3$, $\be =\pi /3$ and
$\al =\pi /6$. That is, the statement of the theorem is true
in this case. Therefore, we may assume that $p>q$ holds for the equation
$p\al +q\be +r\ga =\de _6$ at each vertex of $R_6$.

Then there must be an equation $p\al +q\be +r\ga =\si$ at a vertex different
from the vertices of $R_6$ such that $p<q$. Since $\si =\pi$ or $2\pi$, we get
$(q-p)b+p+r=p(1-b)+qb+r=2$ or $4$. Putting $s=q-p$ we get $sb=u$, where
$u$ equals one of $1,2,3,4$, and $s$ is a positive integer. By $1/2\le b<1$
we obtain $2\le s\le 8$. Since $b=u/s$, the left hand side of \eqref{e2} equals
$v/s$, where $v$ is an integer. Then $v/s=4/3$ gives $s=3$ or $s=6$.

If $s=3$, then $u=2$, $3b=2$, $b=2/3$, $\be =\pi /3$, $\al =\pi /6$ and we are
done. If $s=6$, then $u=3$ or $4$. If $u=3$, then $6b=3$ and $b=1/2$. In this
case, however, the left hand side of \eqref{e2} equals $w/2$, where $w$ is an
integer, which is impossible. If $u=4$, then $6b=4$, $b=2/3$, and we obtain
$\al =\pi /6$ again. This completes the proof. \hfill $\square$

\vfill \eject

\noi
Mikl\'os Laczkovich

\noi
{\tt miklos.laczkovich@gmail.com}

\noi
Department of Analysis, E\"otv\"os Lor\'and University (ELTE), P\'azm\'any
P\'eter s\'et\'any 1/C, 1117 Budapest, Hungary

\noi
Ivan Vasenov

\noi
University of Oxford, Merton College, Merton Street, Oxford, OX1 4JD,
UK

\noi
{\tt vasenovivan@gmail.com}

\end{document}